\documentclass[11pt,a4]{article}
\usepackage{makeidx}
\usepackage{amsmath}
\usepackage{amsfonts,amssymb}
\usepackage{graphicx}

\newtheorem{definition}{Definition}
\newtheorem{theorem}{Theorem}
\newtheorem{lemma}{Lemma}
\newtheorem{proposition}{Proposition}

\input{tcilatex}

\begin{document}

\title{\textbf{On Zero Controllability of Evolution Equations}}
\author{B. Shklyar}
\maketitle

\begin{abstract}
The exact controllability to the origin for linear evolution control
equation is considered.The problem is investigated by its transformation to
infinite linear moment problem.

Conditions for the existence of solution for infinite linear moment problem
has been obtained. The obtained results are applied to the zero
controllability for control evolution equations.
\end{abstract}


\section*{Introduction}

\setcounter{equation}{0} \setcounter{theorem}{0} \setcounter{lemma}{0} %
\setcounter{corollary}{0} \setcounter{definition}{0}\setcounter{equation}{0} %
\setcounter{theorem}{0} \setcounter{lemma}{0} \setcounter{corollary}{0}

Let $X$ be a separable complex Hilbert space.

Given sequences $\left\{ c_{n},n=1,2\ldots ,\ \right\} $ and $\left\{
x_{n}\in X,n=1,2\ldots ,\ \right\} $~find necessary and sufficient
conditions for the existence of an element $g\in X$ such that 
\begin{equation*}
c_{n}=\left( x_{n},g\right) ,n=1,2\ldots ,\ .
\end{equation*}

The problem formulated above is called the linear moment problem. It has a
long history and many applications in geometry, physics, mechanics. 


The goal of this paper is to establish necessary and sufficient conditions
of exact null-controllability for linear evolution control equations with
unbounded input operator by transformation of exact null-controllability
problem (controllability to the origin) to linear infinite moment problem.

It is well-known, that if the sequence $\left\{ x_{n},n=1,2,...,\right\} $%
forms a Riesz basic in the closure of its linear span, the linear moment
problem has a solution if and only if $\sum_{n=1}^{\infty }\left\vert
c_{n}\right\vert ^{2}<\infty $ and vice-versa \cite{Bari}, \cite%
{Gohberg&Krein}, \cite{Ullrich},~\cite{Young80}. This well-known fact is one
of main tools for the controllability analysis of various partial hyperbolic
control equations and functional differential control systems of neutral
type.

However the sequence $\left\{ x_{n},n=1,2,...,\right\} $ doesn't need to be
a Riesz basic for the solvability of linear moment problem. This case
appears under the investigation of the controllability of parabolic control
equations or hereditary functional differential control systems. In this
paper we consider the zero controllability of control evolution equations
for the case when the sequence $\left\{ x_{n},n=1,2,...,\right\} $ of the
moment problem obtained by the transformation of the source control problem
doesn't form a Riesz basic in its closed linear span.

\section{Problem statement}

\setcounter{equation}{0} \setcounter{theorem}{0} \setcounter{lemma}{0} %
\setcounter{corollary}{0} \setcounter{definition}{0}\setcounter{equation}{0} %
\setcounter{theorem}{0} \setcounter{lemma}{0} \setcounter{corollary}{0}

Let $X,U$ be complex Hilbert spaces, and let $\,A$ be infinitesimal
generator of strongly continuous $C_{0}$-semigroups $S\left( t\right) $ in $%
X $ \cite{Hille&Phillips},\cite{Krein}. Consider the abstract evolution
control equation \cite{Hille&Phillips}, \cite{Krein}

\begin{equation}
\dot{x}\left( t\right) =Ax\left( t\right) +Bu\left( t\right) ,x\left(
0\right) =x_{0},\;\,0\leq t<+\infty ,  \label{1.1}
\end{equation}%
where $x\left( t\right) ,\;x_{0}\in X,u\left( t\right) ,u_{0}\in
U,\;B:U\rightarrow X$ is a linear possibly unbounded operator, $W\subset
X\subset V$ are Hilbert spaces with continuous dense injections, where $%
W=D\left( A\right) ~$equipped with graphic norm,$~V=W^{\ast }$, the operator 
$B$ is a bounded operator from $U$ to $V$ (see more details in \cite{Salamon}%
,~ \cite{Da Pratto}, \cite{Nagel}, \cite{Weiss}).

It is well-known that \cite{Da Pratto}, \cite{Nagel},\cite{Salamon}, \cite%
{Weiss}), etc. :

$\bullet $ for each $t\geq 0$ the operator $S\left( t\right) $ has an unique
continuous extension $\mathcal{S}\left( t\right) $ on the space $V$ and the
family of operators $\mathcal{S}\left( t\right) :V\rightarrow V$ is the
semigroup in the class $C_{0}$ with respect to the norm of $V$ and the
corresponding infinitesimal generator $\mathcal{A}$ of the semigroup $%
\mathcal{S}\left( t\right) $ is the closed dense extension of the operator $%
A $ on the space $V$ with domain $D\left( \mathcal{A}\right) =X$;

$\bullet $ the sets of eigenvalues and of generalized eigenvectors of
operators $\mathcal{A},\mathcal{A}^{\ast }$ and $A,\,A^{\ast }$ are the same;

$\bullet $ for each $\mu \notin \sigma \left( A\right) \,$ the resolvent
operator $R_{A}\left( \mu \right) \,$ has a unique continuous extension to
the resolvent operator $\mathcal{R}_{A}\left( \mu \right) :V\rightarrow X$;

$\bullet $ a mild solution $x\left( t,x_{0},u\left( \cdot \right) \right) $
of equation (\ref{1.1}) with initial condition $x\left( 0\right) =x_{0}$ is
obtained by the following representation formula

\begin{equation}
x(t,x_{0},u(\cdot ))=S(t)x_{0}+\int\limits_{0}^{t}\mathcal{S}(t-\tau
)Bu(\tau )d\tau ,  \label{1.2}
\end{equation}%
where the integral in (\ref{2.3}) is understood in the Bochner's sense \cite%
{Hille&Phillips}. To assure $x(t,x_{0},u(\cdot ))\in X,~\forall x_{0}\in
X,u(\cdot )\in L_{2}^{\mathrm{loc}}\left[ 0,+\infty \right) ,t\geq 0,$ we
assume that $\int\limits_{0}^{t}\mathcal{S}(t-\tau )Bu(\tau )d\tau \in X$
for any $u(\cdot )\in L_{2}^{\mathrm{loc}}\left[ 0,+\infty \right) ,t\geq 0$ 
\cite{Salamon}, \cite{Weiss}. 

\begin{definition}
\label{D2.1} Equation (\ref{1.1}) is said to be exact null-controllable on $%
\left[ 0,t_{1}\right] $ by controls vanishing after time moment $t_{2}$ if
for each $x_{0}\in X$ there exists a control $u\left( \cdot \right) \in
L_{2}\left( \left[ 0,t_{2}\right] ,U\right) ,u\left( t\right) =0$ \textrm{%
a.e. }on $[t_{2},+\infty )$ such that 
\begin{equation}
x\left( t_{1},x_{0},u\left( \cdot \right) \right) =0.  \label{1.3}
\end{equation}
\end{definition}

\label{Assumptions}

\subsection{The assumptions}

The assumptions on $A$ are listed below.

\begin{enumerate}
\item The operators $A$ has purely point spectrum $\sigma $ with no finite
limit points. Eigenvalues of $A$ have finite multiplicities.

\item There exists $T\geq 0$ such that \thinspace $\,$all mild solutions of
the equation $\dot{x}\left( t\right) =Ax\left( t\right) $ are expanded in a
series of generalized eigenvectors of the operator $A$ converging uniformly
for any $t\in \left[ T_{1},T_{2}\right] ,T<T_{1}<T_{2}.$
\end{enumerate}

\section{Main results}

\setcounter{equation}{0} \setcounter{theorem}{0} \setcounter{lemma}{0} %
\setcounter{corollary}{0} \setcounter{definition}{0}\setcounter{equation}{0} %
\setcounter{theorem}{0} \setcounter{lemma}{0} \setcounter{corollary}{0}

\subsection{One input case}

For the sake of simplicity we consider the following:

\begin{enumerate}
\item The operator $A$ has all the eigenvalues with multiplicity $1$.

\item $U=\mathbb{R}$ (one input case). It means that the possibly unbounded
operator $B:U\rightarrow \mathbb{R}$ is defined by an element $b\in V$, i.e.
equation (\ref{1.1}) can be written in the form%
\begin{equation}
\dot{x}\left( t\right) =Ax\left( t\right) +bu\left( t\right) ,x\left(
0\right) =x_{0},b\in V,\;\,0\leq t<+\infty .  \label{2.1}
\end{equation}

The operator defined by $b\in V$ is bounded if and only if $b\in X.$\qquad
\end{enumerate}

Let the eigenvalues $\lambda _{j}\in \sigma ,j=1,2,\ldots $ of the operator $%
A$ be enumerated in the order of non-decreasing of their absolute values,
and let $\varphi _{j},\psi _{j},j=1,2,\ldots ,$be eigenvectors of the
operator $A$ and the adjoint operator $A^{\ast }$ respectively. It is
well-known, that 
\begin{equation}
(\varphi _{_{k},}\psi _{j})=\delta _{kj},\ j,k=1,2\ldots ,\   \label{2.2}
\end{equation}%
where $\delta _{kj},\ j,k=1,2\ldots $is the Kroneker delta.

Denote:%
\begin{equation}
x_{j}\left( t\right) =\left( x\left( t,x_{0},u\left( \cdot \right) \right)
,\psi _{j}\right) ,~x_{0j}=\left( x_{0},\psi _{j}\right) ,~b_{j}=\left(
b,\psi _{j}\right) ,~j=1,2,....  \label{2.3}
\end{equation}%
All scalar products in (\ref{2.3}) are correctly defined, because $\psi
_{j}\in W,$ $b\in V=W^{\ast }.$

\begin{theorem}
\label{T2.1}For equation (\ref{1.1}) to be exact null-controllable on $\left[
0,t_{1}\right] ,$ $t_{1}>T,$ by controls vanishing after time moment $%
t_{1}-T $, it is necessary and sufficient that the following infinite moment
problem 
\begin{equation}
x_{0j}=-\int_{0}^{t_{1}-T}e^{-\lambda _{j}\tau }b_{j}u\left( \tau \right)
d\tau ,~j=1,2,...  \label{2.4}
\end{equation}%
with respect to $u\left( \cdot \right) \in L_{2}\left[ 0,t_{1}-T\right] $ is
solvable for any $x_{0}\in X$ .
\end{theorem}

\textbf{Proof.} \textbf{Necessity. }Multiplying (\ref{1.1}) by $\psi _{j},$ $%
j=1,2,...,$and using (\ref{2.3}) we obtain%
\begin{eqnarray}
\dot{x}_{j}\left( t\right) &=&\left( Ax\left( t\right) ,\psi _{j}\right)
+b_{j}u\left( t\right) =\left( x\left( t\right) ,A^{\ast }\psi _{j}\right)
+b_{j}u\left( t\right) =  \notag \\
&=&\lambda _{j}x_{j}\left( t\right) +b_{j}u\left( t\right) ,j=1,2,...,.
\label{2.5}
\end{eqnarray}%
Here $x_{j}\left( t\right) ,\dot{x}_{j}\left( t\right) \ $and $%
b_{j},j=1,2,...,$ are well-defined because $\psi _{j}\in W,\ \dot{x}\left(
t\right) ,Ax\left( t\right) ,b\in V=W^{\ast }.$

From (\ref{2.5}) it follows that%
\begin{equation}
x_{j}\left( t\right) =e^{\lambda _{j}t}\left( x_{j0}+\int_{0}^{t}e^{-\lambda
_{j}t}b_{j}u\left( \tau \right) d\tau \right) ,j=1,2,...,.  \label{2.6}
\end{equation}%
In accordance with the definition of exact null-controllability there exists 
$u\left( \cdot \right) \in L_{2}\left( \left[ 0,t_{1}-T\right] ,U\right)
,u\left( t\right) =0$ \textrm{a.e. }on $[t_{1}-T,+\infty )$ such that (\ref%
{1.3}) holds. Using $u\left( t\right) $ and $t_{1}$ in (\ref{2.6}), we
obtain by (\ref{1.3}) and (\ref{2.5}), that

\begin{equation}
x_{j}\left( t_{1}\right) =e^{\lambda _{j}t_{1}}\left(
x_{j0}+\int_{0}^{t_{1}-T}e^{-\lambda _{j}t}b_{j}u\left( \tau \right) d\tau
\right) =0,j=1,2,...,.  \label{2.7}
\end{equation}
Hence we have (\ref{2.4}) to be true. This proves the necessity.

\textbf{Sufficiency. }Let the control $u\left( \cdot \right) \in L_{2}\left( %
\left[ 0,t_{1}-T\right] ,U\right) ,u\left( t\right) =0$ \textrm{a.e. }on $%
[t_{1}-T,+\infty )$ satisfies (\ref{2.4}). It follows from (\ref{2.4}) and (%
\ref{2.7}) that 
\begin{equation}
x_{j}\left( t_{1}-T\right) =\left( x\left( t_{1}-T\right) ,\psi _{j}\right)
=0,j=1,2,....  \label{2.7.1}
\end{equation}

Denote $z\left( t\right) =x\left( t+t_{1}-T\right) ,~t\geq T.$Obviously, $%
z\left( t\right) $ is a mild solution of the equation $\dot{z}\left(
t\right) =Az\left( t\right) $ with initial condition \ $z\left( 0\right)
=x\left( t_{1}-T\right) ~.$By assumption 3 (see the list of assumptions) $%
z\left( t\right) $ $~$is expanded in a series

\begin{equation}
z\left( t\right) =\sum_{j=1}^{\infty }e^{\lambda _{j}t}\left( x\left(
t_{1}-T\right) ,\psi _{j}\right) ,t\geq T,  \label{2.8}
\end{equation}

so by (\ref{2.7.1}) and (\ref{2.8}) we obtain

$z\left( t\right) =x\left( t+t_{1}-T,x_{0},u\left( \cdot \right) \right)
\equiv 0,t\geq T\Leftrightarrow x\left( t,x_{0},u\left( \cdot \right)
\right) \equiv 0,t\geq t_{1}.$

This proves the sufficiency.

\subsection{Solution of moment problem (\protect\ref{2.4})}

The solvability of moment problem (\ref{2.4}) for each $x_{0}\in X\ $%
essentially depends on the properties of eigenvalues $\lambda _{j},$ $%
~j=1,2,...,.$

If the sequence of exponents $\left\{ e^{-\lambda
_{n}t}b_{n},n=1,2,...,\right\} $forms a Riesz basic in $L_{2}\left[ 0,t_{1}-T%
\right] ,$then the moment problem 
\begin{equation}
c_{j}=-\int_{0}^{t_{1}-T}e^{-\lambda _{j}\tau }b_{j}u\left( \tau \right)
d\tau ,~j=1,2,...  \label{2.9}
\end{equation}
is solvable if and only if 
\begin{equation}
\sum_{j=1}^{\infty }\left\vert c_{j}\right\vert ^{2}<\infty  \label{2.10}
\end{equation}%
There are very large number of papers and books devoted to conditions for
sequence of exponents to be a Riesz basic. All these conditions can be used
for sufficient conditions of zero controllability of equation (\ref{1.1}).
They are very useful for the investigation of the zero controllability of
hyperbolic partial control equations and functional differential control
systems of neutral type \cite{Rabah&Sklyar}.

However moment problem (\ref{2.9}) may also be solvable when the sequence $%
\left\{ e^{-\lambda _{n}t}b_{n},n=1,2,...,\right\} $ doesn't form a Riesz
basic in $L_{2}\left[ 0,t_{1}-T\right] .$ Below we will try to find more
extended controllability conditions which are applicable for the case when
the sequence $\left\{ e^{-\lambda _{n}t}b_{n},n=1,2,...,\right\} $ doesn't
form a Riesz basic in $L_{2}\left[ 0,t_{1}-T\right] .$

\begin{definition}
\label{D2.2}The sequence $\left\{ x_{j}\in X,j=1,2,...,\right\} $ is said to
be minimal, if there no element of the sequence belonging to the closure of
the linear span of others. By other words, 
\begin{equation*}
x_{j}\notin \overline{\mathrm{span}}\left\{ x_{k}\in X,k=1,2,...,k\neq
j\right\} .
\end{equation*}
\end{definition}

The investigation of the controllability problem defined above is based on
the following result of Boas \cite{Boas} (see also \cite{Bari} and \cite%
{Young98}).


\textbf{Theorem} \textit{Let }$x_{j}\in X,j=1,2,...,.$\textit{\ The linear
moment problem }%
\begin{equation*}
c_{j}=\left( x_{j},g\right) ,j=1,2,...
\end{equation*}%
\textit{has a solution }$g\in X$\textit{\ for each square summable sequence }%
$\left\{ c_{j},j=1,2,...\right\} $\textit{\ if and only if there exists a
positive constant }$\gamma $\textit{\ such that all the inequalities }%
\begin{equation}
\gamma \sum_{k=1}^{n}\left\vert c_{k}\right\vert ^{2}\leq \left\Vert
\sum_{j=1}^{n}c_{j}x_{j}\right\Vert ^{2},n=1,2,...,.  \label{2.11}
\end{equation}%
\textit{are valid. }


Let $\left\{ x_{j}\in X,j=1,2,...,\right\} $ a sequence of elements of $X$ ,
and let 
\begin{equation*}
G_{n}=\left\{ \left( x\,_{i},x_{j}\right) ,i,j=1,2,...,n\right\}
\end{equation*}%
be the Gram matrix of $n$ first elements $\left\{ x_{1},...,x_{n}\right\} $
of above sequence. Denote by $\gamma _{n}^{\min }$ the minimal eigenvalue of
the $n\times n$-matrix $G_{n}.$Each minimal sequence $\left\{ x_{j}\in
X,j=1,2,...,\right\} $ is linear independent, hence any first $n~$elements $%
\left\{ x_{1},...,x_{n}\right\} ,$ $n=1,2,...,$ of this sequence are linear
independent, so $\gamma _{n}^{\min }>0,$ $\forall n=1,2,...,.~~$It is easily
to show that the sequence $\left\{ \gamma _{n}^{\min },n=1,2,...,\right\} $
decreases , so there exists $\lim\limits_{n\rightarrow \infty }\gamma
_{n}^{\min }\geq 0.$

\begin{definition}
\label{D2.3}The sequence $\left\{ x_{j}\in X,j=1,2,...,\right\} $ is said to
be strongly minimal, if $\gamma ^{\mathrm{\min }}=\lim\limits_{n\rightarrow
\infty }\gamma _{n}^{\min }>0.$
\end{definition}

It is well-known that for Hermitian $n\times n$-matrix $G_{n}=\left\{ \left(
x_{j},x_{k}\right) ,~j,k=1,2,...,n\right\} $

\begin{equation}
\gamma _{n}^{\mathrm{\min }}\sum_{k=1}^{n}\left\vert c_{k}\right\vert
^{2}\leq \sum_{j=1}^{n}\sum_{k=1}^{n}c_{j}\left( x_{j},x_{k}\right) 
\overline{c_{k}},n=1,2,...,.  \label{2.12}
\end{equation}%
From the well-known formula $\sum_{j=1}^{m}\sum_{k=1}^{m}c_{j}\left(
x_{j},x_{k}\right) \overline{c_{k}}=\left\Vert
\sum_{j=1}^{m}c_{j}x_{j}\right\Vert ^{2},$(\ref{2.11}) and the inequality $%
\gamma _{n}^{\mathrm{\min }}\geq $ $\gamma ^{\mathrm{\min }}>0$ it follows
that 
\begin{equation}
\gamma ^{\mathrm{\min }}\sum_{k=1}^{n}\left\vert c_{k}\right\vert ^{2}\leq
\left\Vert \sum_{j=1}^{n}c_{j}x_{j}\right\Vert ^{2}  \label{2.13}
\end{equation}

Hence the above theorem can be reformulated as follows

\begin{theorem}
\label{T2.2} The linear moment problem 
\begin{equation}
c_{j}=\left( x_{j},g\right) ,j=1,2,...  \label{2.14}
\end{equation}%
has a solution $g\in X$ for any sequence $\left\{ c_{n},n=1,2,...\right\} ,$ 
$\sum\limits_{j=1}^{\infty }c_{j}^{2}<\infty $ if and only if the sequence $%
\left\{ x_{n},n=1,2,..,\right\} $is strongly minimal.
\end{theorem}

\section{Solution of the exact null-controllability problem.}

\setcounter{equation}{0} \setcounter{theorem}{0} \setcounter{lemma}{0} %
\setcounter{corollary}{0} \setcounter{definition}{0}\setcounter{equation}{0} %
\setcounter{theorem}{0} \setcounter{lemma}{0} \setcounter{corollary}{0}

%
%
%

\begin{theorem}
\label{T3.1}For equation (\ref{1.1}) to be exact null-controllable on $\left[
0,t_{1}\right] ,$ $t_{1}>T,$ by controls vanishing after time moment $%
t_{1}-T $, it is necessary, that the sequence 
\begin{equation}
\left\{ e^{-\lambda _{j}\tau }b_{j},t\in \left[ 0,t_{1}-T\right]
,~j=1,2,...,\right\}  \label{3.1}
\end{equation}%
is minimal, and sufficient , that:

\begin{itemize}
\item the sequence $\left\{ e^{-\lambda _{j}\tau }b_{j},t\in \left[ 0,t_{1}-T%
\right] ,~j=1,2,...\right\} $is strongly minimal;

\item 
\begin{equation}
\sum_{j=1}^{\infty }\left\vert \left( x_{0},\psi _{j}\right) \right\vert
^{2}<+\infty ,\forall x_{0}\in X.  \label{3.2}
\end{equation}
\end{itemize}
\end{theorem}

\textbf{Proof. Necessity. }If the problem (\ref{2.4}) has a solution for any 
$x_{0}\in X,$then it has a solution for any eigenvector $\varphi
_{k},k=1,2,...,$ of the operator $A,$ so for each $k=1,2,...,$ there exists
a function $u_{k}\left( \cdot \right) \in L_{2}\left[ 0,t_{1}-T\right] $
such that 
\begin{equation}
\left( \varphi _{k},\psi _{j}\right) =-\int_{0}^{t_{1}-T}e^{-\lambda
_{j}\tau }b_{j}u_{k}\left( \tau \right) d\tau ,~~~j=1,2,...,.  \label{3.3}
\end{equation}%
The sequence $\left\{ \varphi _{k},k=1,2,...,\right\} $ of eigenvectors of
the operator $A$ is biorthogonal to the sequence $\left\{ \psi
_{k},k=1,2,...,\right\} $ of eigenvectors of the operator $A^{\ast }.$ Hence
it follows from (\ref{3.3}) and (\ref{2.2}) that

\begin{equation*}
\delta _{jk}=\left( \varphi _{k},\psi _{j}\right)
=-\int_{0}^{t_{1}-T}e^{-\lambda _{j}\tau }b_{j}u_{k}\left( \tau \right)
d\tau ,j=1,2,...,.
\end{equation*}%
i.e. the sequence $\left\{ -u_{k}\left( t\right) ,t\in \left[ 0,t_{1}-T%
\right] ,~k=1,2,...,\right\} $ is biorthogonal to the sequence$~\left\{
e^{-\lambda _{j}t}b_{j},t\in \left[ 0,t_{1}-T\right] ,~j=1,2,...,\right\} .$

It proves the necessity.

\textbf{Sufficiency. }The sufficiency follows immediately from (\ref{3.2})
and Theorem \ref{T2.2}.

It proves the theorem.

\subsection{The case of the strongly minimal sequence of eigenvectors of the
operator $A$.}

Obviously the sequence of eigenvectors of the operator $A$ being considered
is a minimal sequence.

Below we consider the operator $A$ having the strongly minimal sequence of
eigenvectors.

\begin{theorem}
\label{T3.2}Let the sequence $\left\{ \varphi _{j},j=1,2,...\right\} $ of
eigenvectors of the operator $A$ be strongly minimal.

For equation (\ref{1.1}) to be exact null-controllable on $\left[ 0,t_{1}%
\right] ,$ $t_{1}>T,$ by controls vanishing after time moment $t_{1}-T$, it
is necessary, that the sequence $\left\{ e^{-\lambda _{j}\tau }b_{j},t\in %
\left[ 0,t_{1}-T\right] ,~j=1,2,...\right\} $ is minimal, and sufficient,
that $\func{Re}\lambda _{j}\geq \beta $ for some $\beta \in \mathbb{R}$ and
the sequence$\left\{ e^{-\lambda _{j}t}b_{j},t\in \left[ 0,t_{1}-T\right]
,~j=1,2,...\right\} $ is strongly minimal.
\end{theorem}

\textbf{Proof.} The necessity follows from Theorem \ref{3.1}.

\textbf{Sufficiency. }By Assumption 3 of the list of assumptions the series 
\begin{equation}
\sum\limits_{j=1}^{\infty }\left( x_{0},\psi _{j}\right) e^{\lambda
jt}\varphi _{j},\forall t>T  \label{3.4}
\end{equation}%
converges. Since the sequence $\left\{ \varphi _{j},j=1,2,...\right\} $ of
eigenvectors of the operator $A$ is strongly minimal, then on account of
property (\ref{2.9} there exists a number $\alpha $ such that

\begin{eqnarray}
\alpha ^{2}\sum\limits_{j=1}^{n}\left\vert \left( x_{0},\psi _{j}\right)
\right\vert ^{2}e^{2\func{Re}\lambda _{j}t} &\leq
&\sum_{j=1}^{n}\sum_{k=1}^{n}\left( x_{0},\psi _{j}\right) e^{\lambda
jt}\left( \varphi _{j},\varphi _{k}\right) \overline{\left( x_{0},\psi
_{k}\right) }e^{\overline{\lambda _{k}}t},  \label{3.5} \\
\forall x_{0} &\in &X,~\forall n\in \mathbb{N},~\forall t>T.  \notag
\end{eqnarray}%
It follows from (\ref{3.4}) and (\ref{3.5}) that 
\begin{equation}
\sum\limits_{j=1}^{\infty }\left\vert \left( x_{0},\psi _{j}\right)
\right\vert ^{2}e^{2\func{Re}\lambda _{j}t}<+\infty ,\forall x_{0}\in
X,\forall t>T.  \label{3.6}
\end{equation}

As $\func{Re}\lambda _{j}\geq \beta $ for some $\beta \in \mathbb{R},$we
have by (\ref{3.6}) that (\ref{3.2}) holds.

In accordance with Theorem \ref{T3.1} condition (\ref{3.2}) and the strong
minimality of the sequence (\ref{3.1})~imply the exact null-controllability
of equation (\ref{1.1}). It proves the theorem.

\subsection{The case when the eigenvectors of the operator $A$ form a Riesz
basic}

One of the important problems of the operator theory is the case when the
generalized eigenvectors of the operator $A$\ being considered form a Riesz
basic in $X.$\ The problem of expansion into a Riesz basic of eigenvectors
of the operator $A$ is widely investigated in the literature (see, for
example, \cite{Ahiezer&Glazman}, \cite{Gen&Siu}, \cite{Gohberg&Krein}, \cite%
{Naimark} and references therein). Obviously the sequence of these vectors
is strongly minimal. In this case one can set $T=0,$\ so the Theorems \ref%
{T3.1}, \ref{T3.2} and Lemma \ref{L3.1} can be proven with $T=0.$

\begin{theorem}
\label{T3.3}Let the sequence of operator $A$ forms a Riesz basic in $X.$%

For equation (\ref{1.1}) to be exact null-controllable on $\left[ 0,t_{1}
\right], t_{1}>T,$ by controls vanishing after time moment $t_{1}-T$, it
is necessary and sufficient, that the sequence sequence$\left\{ e^{-\lambda
_{j}t}b_{j},t\in \left[ 0,t_{1}-T\right] ,~j=1,2,...\right\} $ is strongly
minimal .
\end{theorem}

\textbf{Proof.} Let $\left\{ c_{j},j=1,2,...,\right\} $ be any complex
sequence satisfying the condition $\sum_{j=1}^{\infty }\left\vert
c_{j}\right\vert ^{2}<\infty .$

Since the sequence $\left\{ \varphi _{j},j=1,2,...,\right\} $ of
eigenvectors of the operator $A$ forms the Riesz basic, there exists a
vector $x_{0}\in X$ such that 
\begin{equation*}
c_{j}=\left( x_{0},\psi _{j}\right) ,j=1,2,...,
\end{equation*}%
so in virtue of Theorem \ref{T2.1} the exact null controllability being
considered in the paper is equivalent to the solvability of the linear
moment problem 
\begin{equation}
c_{j}=\int_{0}^{t_{1}-T}e^{-\lambda _{j}\tau }b_{j}u\left( \tau \right)
d\tau ,~j=1,2,...,  \label{3.7}
\end{equation}%
for any complex sequence $\left\{ c_{j},j=1,2,...,\right\} ~$satisfying the
condition $\sum_{j=1}^{\infty }\left\vert c_{j}\right\vert ^{2}<\infty .$

By above mentioned results of \cite{Boas} and \cite{Bari} the linear moment
problem (\ref{3.7}) is solvable for any complex sequence $\left\{
c_{j},j=1,2,...,\right\} $ satisfying the condition $\sum_{j=1}^{\infty
}\left\vert c_{j}\right\vert ^{2}<\infty \ $if and only if$\ $the sequence $%
\left\{ e^{-\lambda _{j}t}b_{j},t\in \left[ 0,t_{1}-T\right]
,~j=1,2,...\right\} $ is strongly minimal . It proves the theorem.

Obviously, the condition $b_{j}\neq 0,j=1,2,...,$ is the necessary condition
for the solvability of the moment problem (\ref{2.1}).

\begin{lemma}
\label{L3.1}If the sequence 
\begin{equation}
\left\{ e^{-\lambda _{j}t},t\in \left[ 0,t_{1}-T\right] ,~j=1,2,...\right\}
\label{3.8}
\end{equation}%
is strongly minimal and 
\begin{equation}
\inf\limits_{n\in \mathbb{N}}\left\vert b_{n}\right\vert =\beta >0
\label{3.9}
\end{equation}
holds, then the sequence $\left\{ e^{-\lambda _{j}t}b_{j},t\in \left[
0,t_{1}-T\right] ,~j=1,2,...\right\} $is also strongly minimal.
\end{lemma}

\textbf{Proof. }Let the sequence\textbf{\ }$\left\{ e^{-\lambda _{j}t},t\in %
\left[ 0,t_{1}-T\right] ,~j=1,2,...\right\} $ be strongly minimal. From (\ref%
{2.11}) it follows that

\begin{equation}
\alpha \sum_{k=1}^{n}\left\vert c_{k}\right\vert ^{2}\left\vert
b_{j}\right\vert ^{2}\leq \int_{0}^{t_{1}-T}\left\vert
\sum_{j=1}^{n}c_{j}e^{-\lambda _{j}t}b_{j}\right\vert ^{2}dt  \label{3.10}
\end{equation}%
for some positive $\alpha \ $and for every finite sequence $\left\{
c_{1},c_{2},...,c_{n}\right\} .$ By (\ref{3.9}) and (\ref{3.10}) we have

\begin{equation}
\gamma \sum_{k=1}^{n}\left\vert c_{k}\right\vert ^{2}\leq
\int_{0}^{t_{1}-T}\left\vert \sum_{j=1}^{n}c_{j}e^{-\lambda
_{j}t}b_{j}\right\vert ^{2}dt,n=1,2,...,\gamma =\alpha \beta >0.
\label{3.11}
\end{equation}%
where $\gamma =\alpha \beta >0.$ It proves the lemma.

\textbf{Example of strongly minimal sequence. }Below we will prove that the
sequence $\left\{ e^{n^{2}\pi ^{2}t},n=1,2,...,t\in \left[ 0,t_{1}\right]
\right\} $ is strongly minimal for any $t_{1}>0.$

Let $t_{1}=2t_{2}.$ The series $\sum_{n=1}^{\infty }\frac{1}{n^{2}\pi ^{2}}$
converges and $\left( n+1\right) ^{2}-n^{2}\geq 1$, so the sequence $\left\{
e^{n^{2}\pi ^{2}t},n=1,2,...,t\in \left[ 0,t_{2}\right] \right\} $is minimal 
\cite{Fattorini&Russel}. In virtue of Theorem 1.5 of \cite{Fattorini&Russel}
for each $\varepsilon >0$ there exists a positive constant $K_{\varepsilon }$
such that the biorthogonal sequence $\left\{ w_{n}\left( t\right)
,n=1,2,...,t\in \left[ 0,t_{2}\right] \right\} $ satisfies the condition%
\begin{equation}
\left\Vert w_{n}\left( \cdot \right) \right\Vert <K_{\varepsilon
}e^{\varepsilon n^{2}\pi ^{2}},n=1,2,...,.  \label{3.12}
\end{equation}%
The positive constant $\varepsilon $ can be chosen such that $%
t_{2}-\varepsilon >0.$

By the Minkowsky inequality and (\ref{3.12}) one can show that

$\sum_{n=1}^{p}\sum_{m=1}^{p}c_{n}e^{-n^{2}\pi ^{2}t_{2}}\left(
\int_{0}^{t_{2}}w_{n}\left( t\right) w_{m}\left( t\right) dt\right)
e^{-m^{2}\pi ^{2}t_{2}}c_{m}=\int_{0}^{t_{2}}\left(
\sum_{n=1}^{p}c_{n}e^{-n^{2}\pi ^{2}t_{2}}w_{n}\left( t\right) dt\right)
^{2}dt\leq $

$\leq \int_{0}^{t_{2}}\sum_{n=1}^{p}\left\vert c_{n}\right\vert
^{2}\sum_{n=1}^{p}\left\vert e^{-n^{2}\pi ^{2}t_{2}}w_{n}\left( t\right)
\right\vert ^{2}dt=\sum_{n=1}^{p}\left\vert c_{n}\right\vert
^{2}\sum_{n=1}^{p}e^{-2n^{2}\pi ^{2}t_{2}}\int_{0}^{t_{2}}\left\vert
w_{n}\left( t\right) \right\vert ^{2}dt\leq $

$\leq \sum_{n=1}^{p}\left\vert c_{n}\right\vert
^{2}\sum_{n=1}^{p}e^{-2n^{2}\pi ^{2}t_{2}}\left\Vert w_{n}\left( \cdot
\right) \right\Vert ^{2}\leq K_{\varepsilon }^{2}\sum_{n=1}^{p}\left\vert
c_{n}\right\vert ^{2}\sum_{n=1}^{p}e^{-2n^{2}\pi ^{2}\left(
t_{2}-\varepsilon \right) }.$

The series $\sum_{n=1}^{\infty }e^{-2n^{2}\pi ^{2}\left( t_{2}-\varepsilon
\right) }$ converges for any $t_{2},\varepsilon ,t_{2}>\varepsilon ,~$so $%
\sum_{n=1}^{p}e^{-2n^{2}\pi ^{2}\left( t_{2}-\varepsilon \right) }\leq M$,
where $M$ is a positive constant.

Hence

\begin{equation}
\sum_{n=1}^{p}\sum_{m=1}^{p}c_{n}e^{-n^{2}\pi ^{2}t_{2}}\left(
\int_{0}^{t_{2}}w_{n}\left( t\right) w_{m}\left( t\right) dt\right)
e^{-m^{2}\pi ^{2}t_{2}}c_{m}\leq K_{\varepsilon
}^{2}M\sum_{n=1}^{p}\left\vert c_{n}\right\vert ^{2}  \label{3.13}
\end{equation}%
for every finite sequence $\left\{ c_{1},c_{2},...,c_{p}\right\} .$
Obviously the sequence

$\left\{ h_{n}\left( t\right) =\left\{ 
\begin{array}{cc}
e^{-n^{2}\pi ^{2}t_{2}}w_{n}\left( t-t_{2}\right) , & t\in \left[
t_{2},2t_{2}\right] , \\ 
0, & t\in \left[ 0,t_{2}\right)%
\end{array}%
\right. ,n=1,2,...,\right\} $is the biorthogonal sequence to the sequence $%
\left\{ e^{n^{2}\pi ^{2}t},n=1,2,...,t\in \left[ 0,t_{1}\right] \right\} ,$
and

$\left( \int_{0}^{t_{1}}h_{n}\left( t\right) h_{m}\left( t\right) dt\right)
=e^{-n^{2}\pi ^{2}t_{2}}\left( \int_{t_{2}}^{2t_{2}}w_{n}\left(
t-t_{2}\right) w_{m}\left( t-t_{2}\right) dt\right) e^{-m^{2}\pi
^{2}t_{2}}=e^{-n^{2}\pi ^{2}t_{2}}\left( \int_{0}^{t_{2}}w_{n}\left(
t\right) w_{m}\left( t\right) dt\right) e^{-m^{2}\pi ^{2}t_{2}},$

\noindent so it follows from (\ref{3.13}) that 
\begin{equation*}
\sum_{n=1}^{p}\sum_{m=1}^{p}c_{n}\left( \int_{0}^{t_{1}}h_{n}\left( t\right)
h_{m}\left( t\right) dt\right) c_{m}\leq K_{\varepsilon
}^{2}M\sum_{n=1}^{p}\left\vert c_{n}\right\vert ^{2}.
\end{equation*}%
Hence \cite{Kaczmarz&Steinhaus} 
\begin{equation}
\sum_{n=1}^{p}\sum_{m=1}^{p}c_{n}\left( \int_{0}^{2t_{1}}e^{n^{2}\pi
^{2}\tau }e^{m^{2}\pi ^{2}\tau }\right) c_{m}d\tau \geq \gamma
\sum_{n=1}^{p}\left\vert c_{n}\right\vert ^{2},p=1,2,...,  \label{3.14}
\end{equation}%
for every finite sequence $\left\{ c_{1},c_{2},...,c_{p}\right\} ,$where $%
\gamma =\frac{1}{K_{\varepsilon }^{2}M}>0.$ It proves that the sequence $%
\left\{ e^{n^{2}\pi ^{2}t},t\in \left[ 0,t_{1}\right] ,~n=1,2,...\right\} $
is strongly minimal for any $~t_{1}>0$.\textit{\ }

\section{Approximation Theorems}

\setcounter{equation}{0} \setcounter{theorem}{0} \setcounter{lemma}{0} %
\setcounter{corollary}{0} \setcounter{definition}{0}\setcounter{equation}{0} %
\setcounter{theorem}{0} \setcounter{lemma}{0} \setcounter{corollary}{0}

As was said at the end of the previous section the condition $%
\lim\limits_{n\rightarrow \infty }\lambda _{n}^{\min }$ $>0$ in general can
be checked by numerical methods. The problem appears to be rather difficult
in general.

However there are sequences for which the validity of above inequality can
be easily established. For example, every orthonormal sequence is strongly
minimal.

Below we will show that if the sequence 
\begin{equation*}
\left\{ y_{j}\in X,j=1,2,...\right\}
\end{equation*}%
can be approximated in the some sense by strongly minimal sequence 
\begin{equation*}
\left\{ x_{j}\in X,j=1,2,...\right\} ,
\end{equation*}%
then it is also strongly minimal.

\begin{theorem}
\label{T4.1}If the sequence $\left\{ x_{j}\in X,j=1,2,...\right\} $ is
strongly minimal, let the sequence $\left\{ y_{j}\in X,j=1,2,...\right\} $
be such that the sequence $\left\{ P_{n}y_{j}-x_{j},j=1,2,...\right\} $ is
linear independent and 
\begin{equation}
\left\Vert \sum_{j=1}^{n}c_{j}\left( y_{j}-x_{j}\right) \right\Vert \leq
q\left\Vert \sum_{j=1}^{n}c_{j}x_{j},\right\Vert ,n=1,2,...~,  \label{4.1}
\end{equation}%
where $\left\{ c_{j},j=1,2,...\right\} $ is any sequence of complex numbers, 
$q$ is a constant, $0<q<1,$ then the sequence $\left\{ y_{j}\in
X,j=1,2,...\right\} $ also is strongly minimal.
\end{theorem}

\textbf{Proof}. 
%
%
Let $\left\{ c_{k},k=1,2,...\right\} $ be an arbitrary sequence of complex
number. Denote:

\begin{equation}
x^{0}=\sum_{k=1}^{n}c_{k}x_{k},~x^{1}=\sum_{k=1}^{n}c_{k}\left(
x_{k}-y_{k}\right) .  \label{4.2}
\end{equation}

From (\ref{4.2}) it follows, that

\begin{equation}
x^{0}=x^{1}+\sum_{k=1}^{n}c_{k}y_{k},~n=\ 1,2,....  \label{4.3}
\end{equation}

By (\ref{4.1}) $\ $we obtain that

\begin{equation}
\left\Vert x^{1}\right\Vert \leq q\left\Vert x^{0}\right\Vert .  \label{4.4}
\end{equation}

Hence using (\ref{4.4}) in (\ref{4.3}) we obtain%
\begin{equation}
\left\Vert x^{0}\right\Vert \leq \frac{1}{1-q}\left\Vert
\sum_{k=1}^{n}c_{k}y_{k}\right\Vert ,~n=\ 1,2,....  \label{4.5}
\end{equation}

Since the sequence $\left\{ x_{j}\in X,j=1,2,...\right\} $ is strongly
minimal and $x^{0}$ $=\sum_{k=1}^{n}c_{k}x_{k}$, we have 
\begin{equation}
\sum_{k=1}^{n}\left\vert c_{k}\right\vert ^{2}\leq \frac{1}{\alpha ^{2}}%
\left\Vert x^{0}\right\Vert ^{2},~n=1,2,...,  \label{4.6}
\end{equation}%
for some $\alpha >0.$

By (\ref{4.6}) and (\ref{4.5}) we obtain $\alpha
^{2}\sum_{k=1}^{n}\left\vert c_{k}\right\vert ^{2}\leq \frac{1}{1-q}%
\left\Vert \sum_{k=1}^{n}c_{k}y_{k}\right\Vert ,~n=\ 1,2,...,$ so

\begin{equation}
\alpha ^{2}\left( 1-q\right) ^{2}\left( \sum_{k=1}^{n}\left\vert
c_{k}\right\vert ^{2}\right) \leq \left\Vert
\sum_{k=1}^{n}c_{k}y_{k}\right\Vert ,~n=\ 1,2,...,.  \label{4.7}
\end{equation}%
Using in (\ref{4.7}) the formula (\ref{2.13}) we obtain

\begin{equation}
\gamma \left( \sum_{k=1}^{n}\left\vert c_{k}\right\vert ^{2}\right) \leq
\sum_{k=1}^{n}\sum_{l=1}^{n}c_{k}\left( y_{k},y_{l}\right) \overline{c_{l}}%
,\gamma =\alpha ^{2}\left( 1-q\right) ^{2}>0  \label{4.8}
\end{equation}%
Let $\mu _{\min }^{\left[ n\right] }$ be a minimal eigenvalue of the Gram
matrix $G_{n}=\left\{ \left( y_{k},y_{l}\right) ,k,l=1,2....\right\} $ for
the sequence $\left\{ y_{j},j=1,2,...,n\right\} .$ From (\ref{4.8}), it
follows that $\lim_{n\rightarrow \infty }\mu _{\min }^{\left[ n\right] }\geq
\gamma >0.$

This proves the theorem.

\subsection{Example}

Let $X=$ $l_{2}$ be the Hilbert space of square summable sequences. Consider
the evolution system

\begin{equation}
\left\{ 
\begin{array}{ccc}
\dot{x}_{k}\left( t\right) =\lambda _{k}x_{k}\left( t\right) +u\left(
t\right) , & k=1,2,..., & 0<t<t_{1}, \\ 
x_{k}\left( 0\right) =x_{k0},n=1,2,..., & k=1,2,..., & 
\end{array}%
\right.  \label{4.9}
\end{equation}%
where $u\left( t\right) ,0<t<t_{1}$ is a scalar control function,$~$

$\left\{ x_{k}\left( t\right) ,k=1,2,...,\right\} ,\left\{
x_{k0},k=1,2,...,\right\} \in l^{2},$ the complex numbers $\lambda _{k},$ $%
k=1,2,...,$belong to the strip $\left\{ z\in \mathbb{C}:\left\vert \func{Re}%
z\right\vert \leq \gamma \right\} ,$ i.e. $\left\vert \func{Re}\lambda
_{k}\right\vert \leq \gamma ,k=1,2,...,$ .

\begin{definition}
\label{D4.1}Equation (\ref{4.9}) is said to be exact null-controllable on $%
\left[ 0,t_{1}\right] $ by controls vanishing after time moment $t_{2},$ if
for each $x_{0}\left( \cdot \right) =\left\{ x_{k0},k=1,2,...,\right\} \in
l_{2}$ there exists a control $u\left( \cdot \right) \in L_{2}\left[ 0,t_{2}%
\right] ,u\left( t\right) =0$ \textrm{a.e. }on $[t_{2},+\infty )$ such that 
\begin{equation*}
x_{k}\left( t\right) \equiv 0,~k=1,2,...,\forall t\geq t_{1}.
\end{equation*}
\end{definition}

Control problem (\ref{4.9}) can be written in the form of (\ref{1.1}), where 
$x\left( t\right) =\left\{ x_{k}\left( t\right) ,k=1,2,...,\right\} \in
l^{2},u\left( \cdot \right) \in L_{2}\left[ 0,t_{1}\right] $; the
self-adjoint operator $A:l_{2}\rightarrow l_{2}$ is defined for $x=\left\{
x_{k},k=1,2,...,\right\} \in l_{2}~$by 
\begin{equation}
Ax=\left\{ \lambda _{k}x_{k},k=1,2,...,\right\}  \label{4.10}
\end{equation}%
with domain $D\left( A\right) =\left\{ x\in l_{2}:Ax\in l_{2}\right\} $, and
the unbounded operator $B$ is defined by 
\begin{equation}
Bu=bu,u\in \mathbb{R},  \label{4.11}
\end{equation}%
where $b=\{1,1,...,1,...\}\notin l_{2}$.

One can show that all the assumptions imposed on equation (\ref{1.1}) are
fulfilled for equation (\ref{4.9}) with $T=0$.

Obviously, the numbers $\lambda _{k},$ $k=1,2,...,$ are eigenvalues of the
operator $A$ defined above; the sequences $e_{k}=\left\{ \underset{1~\mathrm{%
on~}k\text{-\textrm{th place}}}{\underbrace{0,...,0,1,0,...,0}}\right\} $
are corresponding eigenvectors, forming the Riesz basic of $l_{2},$ so $%
b_{j}=1,j=1,2,...,.$

Together with system (\ref{4.9}) consider the other evolution system 
\begin{equation}
\left\{ 
\begin{array}{ccc}
\dot{x}_{k}\left( t\right) =\mu _{k}x_{k}\left( t\right) +u\left( t\right) ,
& n=1,2,..., & 0<t<t_{1}, \\ 
x_{k}\left( 0\right) =x_{k0},k=1,2,..., & n=1,2,..., & 
\end{array}%
\right.  \label{4.12}
\end{equation}%
where 
\begin{equation}
\mu _{k}=\lambda _{k}+O\left( \frac{1}{k}\right) ,k=1,2,...,.  \label{4.13}
\end{equation}%
%
%
%
%
%
%

\begin{proposition}
If system (\ref{4.9}) is exact null-controllable on $\left[ 0,t_{1}\right] $
by controls vanishing after time moment $t_{2},$then the same is valid for
system (\ref{4.12}).
\end{proposition}

\textbf{Proof. }From the Caushy-Schvartz inequality it follows that

$\int_{0}^{t_{2}}\left\vert \sum_{k=1}^{n}c_{k}\left( e^{-\mu
_{k}t}-e^{-\lambda _{k}t}\right) \right\vert ^{2}dt\leq
\sum_{k=1}^{n}\left\vert c_{k}\right\vert
^{2}\int_{0}^{t_{2}}\sum_{k=1}^{n}\left\vert e^{-\mu _{k}t}-e^{-\lambda
_{k}t}\right\vert ^{2}dt=$

$=\sum_{k=1}^{n}\left\vert c_{k}\right\vert
^{2}\int_{0}^{t_{2}}\sum_{k=1}^{n}e^{-2\lambda _{k}t}\left\vert e^{O\left( 
\frac{1}{k}\right) t}-1\right\vert ^{2}dt\leq \sum_{k=1}^{n}\left\vert
c_{k}\right\vert ^{2}\int_{0}^{t_{2}}e^{2\gamma t}\sum_{k=1}^{n}\left\vert
e^{O\left( \frac{1}{k}\right) t}-1\right\vert ^{2}dt.$

The series $\sum_{k=1}^{\infty }\left\vert e^{O\left( \frac{1}{k}\right)
t}-1\right\vert ^{2}$ converges for any $t\geq 0$. Denote 
\begin{equation}
M\left( t_{2}\right) =\int_{0}^{t_{2}}e^{2\gamma t}\sum_{k=1}^{\infty
}\left\vert e^{O\left( \frac{1}{k}\right) t}-1\right\vert ^{2}dt.
\label{4.14}
\end{equation}%
Hence

\begin{equation}
\int_{0}^{t_{2}}\left\vert \sum_{k=1}^{n}c_{k}\left( e^{-\mu
_{k}t}-e^{-\lambda _{k}t}\right) \right\vert ^{2}dt\leq M\left( t_{2}\right)
\sum_{k=1}^{n}\left\vert c_{k}\right\vert ^{2}.  \label{4.15}
\end{equation}%
By Theorem \ref{T3.2} we have the sequence$\left\{ e^{-\lambda _{j}t},t\in %
\left[ 0,t_{2}\right] ,~j=1,2,...\right\} $ to be strongly minimal, so

\begin{equation}
\sum_{k=1}^{n}\left\vert c_{k}\right\vert ^{2}\leq \frac{1}{\alpha ^{2}}%
\int_{0}^{t_{2}}\left\vert \sum_{k=1}^{n}c_{k}e^{-\lambda _{k}t}\right\vert
^{2}dt~\mathrm{for~some~}\alpha >0.  \label{4.16}
\end{equation}

Joining (\ref{4.15}) and (\ref{4.16}) we obtain%
\begin{equation}
\int_{0}^{t_{2}}\left\vert \sum_{k=1}^{n}c_{k}\left( e^{-\mu
_{k}t}-e^{-\lambda _{k}t}\right) \right\vert ^{2}dt\leq
q\int_{0}^{t_{2}}\left\vert \sum_{k=1}^{n}c_{k}e^{-\lambda _{k}t}\right\vert
^{2}dt,~  \label{4.17}
\end{equation}%
where~$q=\frac{M\left( t_{2}\right) }{\alpha }.$

Since from (\ref{4.14}) it follows that $\lim\limits_{t_{1}\rightarrow
\infty }M\left( t_{2}\right) =0,$ one can choose the number $t_{2}$ such
that $0<q<1.$ Hence conditions (\ref{4.17}) are the same as (\ref{4.1}) for $%
x_{k}=e^{-\lambda _{k}t},y_{k}=e^{-\mu _{k}t},k=1,2,...,t\in \left[ 0,t_{2}%
\right] ;q=\frac{M\left( t_{2}\right) }{\alpha ^{2}}.$

As it was said abov by Theorem \ref{T3.2} we have the sequence$\left\{
e^{-\lambda _{j}t},t\in \left[ 0,t_{2}\right] ,~j=1,2,...\right\} $ to be
strongly minimal .

In accordance with Theorem \ref{T4.1} the sequence $\left\{ y_{k}=e^{-\mu
_{k}t},k=1,2,...,t\in \left[ 0,t_{2}\right] \right\} $ is also strongly
minimal, provided that $t_{2}$ is chosen such that $\frac{M\left(
t_{2}\right) }{\alpha ^{2}}<1.$ In accordance with Theorem \ref{T3.1} the
strong minimality of the sequence $\left\{ y_{k}=e^{-\mu
_{k}t},k=1,2,...,t\in \left[ 0,t_{2}\right] \right\} $ provides the zero
controllability of equation (\ref{4.11}) on $\left[ 0,t_{1}\right] $ by
controls vanishing after time moment $t_{2},~\frac{M\left( t_{2}\right) }{%
\alpha ^{2}}<1,$ for any $t_{1}\geq t_{2}$.

\label{References}


\begin{thebibliography}{99}
\bibitem{Ahiezer&Glazman} N. Ahiezer, I. Glazman, \emph{Linear Operator
Theory in Hilbert Spaces}, Moscow, Nauka Publisher, 1966 (Russian).



\bibitem{Boas} R. Boas, A general moment problem, \textsl{Amer. J. Math.}, 
\textbf{63}(1941), 361---370.

\bibitem{Bari} N. Bari, Biorthogonal sequences and bases in Hilbert spaces. 
\textsl{Uchen. Zap. Mosk. Univ.}, 148, Nat, \textbf{4}(1951), 69---107.



\bibitem{Da Pratto} Da Pratto, \textsl{Abstract differential equations and
extrapolation spaces}, Lecture Notes in Mathematics, 1184, Springer-Berlag,
Berlin, New York, 1984.




\bibitem{Fattorini&Russel} H. Fattorini, D. Russel, Uniform bounds on
biorthogonal functions for real exponents with an application to the control
theory of parabolic equations, \textsl{Quart. Appl.Math., }1074, 45 --- 69.

\bibitem{Gen&Siu} Gen Qi Xu, Siu Pang Yung, \textsl{The expansion of
semogroup and a Riesz basic criterion}, J. Diff. Eqn., \textbf{210}(2005), 1
--- 24.

\bibitem{Gohberg&Krein} I. Gohberg, M. Krein,\textsl{\ Introduction to the
Theory of Linear Nonselfadjoint operators, }Transl. math. Monogr., 18, AMS,
Providence, RI, 1969.

\bibitem{Hille&Phillips} E. Hille, R. Philips, \textsl{Functional Analysis
and Semi-Groups}, AMS, 1957.






\bibitem{Kaczmarz&Steinhaus} S. Kaczmarz, H. Steinhaus, \textsl{Theory of
orthogonal series} Monographs Mat., Bd. 6, (PWN, Warsaw), 1958

\bibitem{Krein} M. Krein, \textsl{Linear Differential Equations in Banach
Spaces}, Moscow, Nauka Publisher, 1967 (in Russian).





\bibitem{Nagel} R. Nagel, \textsl{One-parameter semigroups of positive
operators,} Lecture Notes in Notes in Mathematics, 1184, Springer-Berlag,
Berlin, New York, 1984.

\bibitem{Naimark} M. Naimark, \textsl{Linear differential Operators, }%
Moscow, Nauka Publisher, 1969 (in Russian).


\bibitem{Rabah&Sklyar} R. Rabah, G. Sklyar, Thw analysis od exact
controllability of neutral-type systems by the moment problem approach, 
\textsl{SIAM J. Contr. Optimiz.}, \textbf{36} (2007), 2148 --- 2181.


\bibitem{Salamon} D. Salamon, Infinite dimensional linear systems with
unbounded control and observation: a functional analytic approach, \textsl{%
Trans. Amer. Math. Soc.}, \textbf{300}(1987), 383 --- 431.



\bibitem{Weiss} G. Weiss, Admissibility of unbounded control operators, 
\textsl{SIAM J. Contr. and Optimiz}., \textbf{27}(1989), 527 --- 545.

\bibitem{Ullrich} D. Ullrich, Divided differences and systems of nonharmonic
Fourier series, \textsl{Proc. Amer. Math. Soc.}, \textbf{80}(1980), 47 ---
57.

\bibitem{Young80} R. Young, \textsl{An Introduction to Nonharmonic Analysis, 
}Academic Press, New York, 1980.

\bibitem{Young98} R. Young, On a class of Riesz-Fisher sequences, \textsl{%
Proceedings of AMS, }\textbf{126}(1998), 1139---1142.
\end{thebibliography}
\end{document}